\documentstyle{article}
\def\Bbb{\bf}
\def\cp{{\Bbb CP}}
\def\bcp{\cp}
\def\kod{\mbox{Kod}}
\def\bea{\begin{eqnarray*}}
\def\eea{\end{eqnarray*}}

\newtheorem{main}{Theorem}

\newtheorem{defn}{Definition}
\newtheorem{thm}{Theorem}
\newtheorem{prop}{Proposition}
\newtheorem{lem}{Lemma}
\newtheorem{cor}{Corollary}
\newtheorem{conj}{Conjecture}
\newenvironment{proof}{\medskip {\bf Proof.}}{\hfill \rule{.5em}{1em}
\\}
\def\qed{\hfill \rule{.5em}{1em}
\\}
\begin{document}
\sloppy
\title{Kodaira Dimension and the Yamabe Problem}
\author{Claude LeBrun\thanks{Supported 
in part by  NSF grant DMS-9505744.} 
\\ 
SUNY Stony
 Brook 
  }
\date{}
\maketitle

\begin{abstract}  The {\em Yamabe invariant} $Y(M)$
of a 
smooth compact manifold is roughly  the supremum of 
the scalar curvatures of unit-volume constant-scalar curvature
Riemannian metrics $g$ on $M$. (To be absolutely
precise, one only considers constant-scalar-curvature 
metrics which are {\em Yamabe minimizers}, but this 
does not affect the {\sl sign} of the answer.)
If $M$ is the underlying smooth 4-manifold of a
  complex algebraic surface $(M,J)$,
it is shown that the sign of $Y(M)$ is 
completely determined by the Kodaira dimension 
  $\kod (M,J)$. More precisely, 
$Y(M) <  0$ iff    $\kod (M,J)=2$; 
$Y(M) = 0$ iff   $\kod (M,J)=0$ or $1$;
and $Y(M) > 0$ iff   $\kod (M,J)=-\infty$. 
 \end{abstract}

\section{Introduction}
 One may define an interesting and natural diffeomorphism invariant
of a compact smooth $n$-manifold $M$  as the supremum
of the scalar curvatures of unit-volume
 constant-scalar-curvature
metrics on $M$. A minor refinement of this definition,
which does not change the sign of the invariant but guarantees
that it is finite, restricts the class of allowed 
constant-scalar-curvature metrics to the so-called
Yamabe minimizers.  We will refer to the resulting
invariant as the {\em Yamabe invariant} $Y(M)$ of
our manifold. For a more precise definition, 
see \S \ref{yam} below.

For 2-manifolds, this invariant is easy to compute;
 the classical Gauss-Bonnet theorem tells us
that the Yamabe invariant of  a smooth compact surface is
just $4\pi \chi$, where $\chi$ is the  
Euler characteristic.
Many important complex-analytic  invariants of 
compact complex curves are  thus determined by the Yamabe invariant 
of the underlying smooth 2-manifold.

One of the most important invariants of 
a compact complex manifold $(M,J)$ is its
{\em Kodaira dimension} $\kod (M)$. 
If $M$ has complex dimension $m$, recall that
the {\em canonical line bundle} $K \to M$ 
is   $K=\wedge^{m,0}$, so that the
holomorphic  sections of $K$ are exactly the 
holomorphic $m$-forms on $(M,J)$. 
One now defines the Kodaira dimension as 
$$\kod (M) = \limsup \frac{\log h^0 (M,K^{\otimes \ell})}{\log \ell}.$$
This can be shown to coincide with the 
maximal complex dimension of the image of $M$
under   pluri-canonical maps to complex projective space,
so that $\kod (M)\in \{ -\infty , 0, 1, \ldots, m\}$. 
A compact complex $m$-manifold is said to be of {\em general type}
if $\kod M =m$. 

For  Riemann surfaces,
the trichotomy $\kod (M) = -\infty , 0, 1$, exactly
coincides with that given by $Y(M) > 0, = 0, < 0$.
The purpose of this paper is to show that
much the same thing happens for compact complex surfaces
$(M^4,J)$. Our main result is 

\begin{main} \label{main} 
Let $M$ be the underlying 4-manifold of a
 compact complex  surface $(M^4,J)$ with $b_1(M)$ even.  
Then 
$$
Y(M) = \left\{
\begin{array}{ll} Y(M) > 0 &\mbox{iff } \kod (M,J) = -\infty \\
Y(M) = 0 &\mbox{iff } \kod (M,J) = 0 \mbox{ or } 1 \\
  Y(M)<  0 &\mbox{iff } \kod (M,J) = 2.
\end{array}\right.
$$ 
\end{main}

The hypothesis that $b_1(M)\equiv 0\bmod 2$
is equivalent to requiring that $(M,J)$ admit 
a K\"ahler metric. It is also equivalent to 
requiring that $(M,J)$ be a deformation of a
complex algebraic surface. Many of the results
contained in this article will also apply in the
non-K\"ahler case, but the overall picture 
remains considerably less clear when $b_1$ is odd. 

One of the main tools in the proof of the above
result is Seiberg-Witten theory. One of the
 most striking consequences of this
theory is that any two diffeomorphic complex
algebraic surfaces must have the same Kodaira
dimension \cite{FM,FQ}.  Since the Yamabe invariant is
obviously a diffeomorphism invariant, 
Theorem  \ref{main}  casts
this result in an interesting new light. 
On the other hand, the Yamabe invariant 
does not by itself distinguish Kodaira
dimension $0$ from Kodaira dimension $1$.
The distinction between these cases instead only
emerges  when 
one asks a finer question:  when is the Yamabe invariant 
actually the scalar curvature of some unit-volume
metric?

 Finally, it should be emphasized that,
despite the cited analogy between complex 
dimensions $1$ and $2$,  the phenomenon 
 explored here does not persist in higher dimensions. 
For example, consider a non-singular complex hypersurface of 
 degree $m+3$ in $\cp_{m+1}$, where $m\geq 3$. 
This complex $m$-manifold has ample 
 canonical bundle, so its Kodaira dimension is $m$. 
However, the underlying smooth manifold is 
simply connected, non-spin, and has real dimension 
$> 4$,   so a surgical  construction of 
Gromov and Lawson \cite{GL} shows that it 
admits metrics of positive scalar curvature,
and thus has $Y >0$. 
Thus  the Yamabe invariant and Kodaira dimension
are correlated only in complex dimensions $1$ and $2$.

\section{Yamabe Invariants} 
\label{yam}

This section is expository in nature,
and is intended as a convenient summary of  
the Yamabe folklore which will be needed in the
rest of the paper. 

Let $M$ be a smooth compact manifold. 
A {\em conformal class} on $M$
is by definition a collection of 
smooth Riemannian metrics on $M$ of the form 
$$[g]=\{ vg ~|~v: M\to {\Bbb R}^+\},$$ 
where $g$ is some fixed Riemannain metric. 
To each such conformal class, one can 
associate a number $Y_{[g]}$, called the
{\em Yamabe constant} of the class, by 
$$Y_{[g]} = \inf_{{\hat{g}}\in [g] } \frac{\int_M 
s_{\hat{g}}~d\mu_{\hat{g}}}{\left(\int_M 
d\mu_{\hat{g}}\right)^{\frac{n-2}{n}}},$$
where $s_{\hat{g}}$ and $d\mu_{\hat{g}}$ respectively
 denote the scalar curvature and volume measure of
${\hat{g}}$.  If $g$ is a Riemannian metric 
for which the scalar curvature is everywhere positive,
or everywhere zero, or everywhere negative, 
one can show that $Y_{[g]}$ is correspondingly 
positive, zero, or negative.

While $Y_{[g]}$ has been defined as 
an infimum, we could have instead defined it
as a {\sl minimum}, for a remarkable theorem
\cite{aubin,lp,rick} 
of Yamabe, Trudinger, Aubin, and Schoen 
asserts that any conformal class 
$[g]$ contains a 
metric  
which actually minimizes the relevant functional.  
Such a metric will be called a {\em Yamabe minimizer}.
Any Yamabe minimizer has constant scalar curvature;
conversely, a constant-scalar-curvature metric $g$ 
with $s_g \leq 0$ is 
automatically a Yamabe minimizer.
 Aubin's piece of  
the proof 
hinges on the observation
that any conformal class on any 
$n$-manifold automatically
satisfies
$Y_{[g]}\leq n(n-1) (V_n)^{2/n}$, where 
$V_n$ is the  volume of 
the  unit $n$-sphere
$S^n\subset {\Bbb R}^{n+1}$.

Given any smooth compact $n$-manifold, we can therefore
\cite{okob,sch}
define an associated real number by 
$$Y(M) = \sup_{[g]} Y_{[g]} = 
\sup_{[g]}\inf_{{\hat{g}}\in [g] } \frac{\int_M
s_{\hat{g}}~d\mu_{\hat{g}}}{\left(\int_M 
d\mu_{\hat{g}}\right)^{\frac{n-2}{n}}}.$$
By construction, this is a diffeomorphism 
invariant of $M$, and it will be called the 
  {\em Yamabe invariant} of $M$ in this paper. 
Notice that $Y(M) > 0$ iff $M$ admits a 
metric of positive scalar curvature. 
If $Y(M) \leq 0$, the invariant
is simply the   supremum of 
the scalar curvatures of unit-volume 
constant-scalar-curvature metrics on $M$, 
since any constant-scalar-curvature metric
of non-positive scalar curvature is automatically a 
Yamabe minimizer. 

The minimax definition of $Y(M)$ is rather unwieldy for
many purposes.   Fortunately,
it can often be calculated in a simpler manner.

\begin{lem}
Let $[g]$ be any conformal class on a smooth compact 
$n$-manifold $M$. Then
$$\inf_{{\hat{g}}\in [g]} \int_M |s_{\hat{g}}|^{n/2} d\mu_{\hat{g}} 
= |Y_{[g]}|^{n/2}.$$
\label{yobo}
\end{lem}
\begin{proof} If $n=2$, the Gauss-Bonnet theorem
tells us that $\int |s|d\mu \geq 
|\int s~d\mu |= |4\pi \chi (M)| = 
|Y_{[g]}|$, with equality iff $s$ does not 
change sign. We  may therefore assume henceforth that
$n \geq 3$. We will break this into two cases, depending
on the sign of $Y_{[g]}$. 

First suppose that $Y_{[g]} \geq 0$ and $n \geq 3$. The H\"older
inequality tells us that  each metric ${\hat{g}}$ satisfies 
$$\left(\int_M |s_{\hat{g}}|^{n/2} d\mu_{\hat{g}} \right)^{2/n}\geq 
\frac{\int_M s_{\hat{g}}d\mu_{\hat{g}}}{\left(\int_M 
 d\mu_{\hat{g}} \right)^{1-\frac{2}{n}}},$$
with equality iff $s_{\hat{g}}$ is a non-negative constant. Since the 
infimum of the right-hand side over ${\hat{g}}\in [g]$ is $Y_{[g]}$
and is achieved by a metric of constant scalar curvature $s\geq 0$, 
the claim follows. 

Finally, we come to the   case in which  $Y_{[g]} < 0$
and $n \geq 3$. 
Let $g\in [g]$ be a metric of constant negative scalar 
curvature, and express any other metric in $[g]$
as  $\hat{g}=
u^\ell g$,
where $\ell = 4/(n-2)$ and $u$ is some smooth
positive function. The scalar curvatures $s=s_g$ and 
$\hat{s}=s_{\hat{g}}$ are related by 
$$\hat{s}u^{\ell + 1}= su + (n-1)\ell \triangle u$$
where $\triangle=d^\ast d$ is the positive Laplacian of $g$. 
Thus
\bea
\int_M \hat{s}u^\ell d\mu  &=&
\int_M  \left(s+(n-1)\ell \frac{1}{u}d^\ast du\right) d\mu\\
&=& \int_M  \left(s-(n-1)\ell 
\frac{|d u |^2}{u^2}\right) d\mu\\
&\leq& \int_Ms~d\mu ,
\eea 
where $d\mu=d\mu_g$. 
The H\"older inequality therefore implies 
 \bea
\left(\int_M |\hat{s}u^\ell |^{n/2}d\mu\right)^{2/n}
\left(\int_Md\mu\right)^{\frac{n-2}{n}}
 &\geq&  \int (-\hat{s}u^\ell ) d\mu\\
&\geq & -\int s ~d\mu\\
&=&|Y_{[g]}|
\left(\int_Md\mu\right)^{(n-2)/n}.
\eea
Hence  
$$
\int_M |s_{\hat{g}}|^{n/2}d\mu_{\hat{g}}=
\int_M |\hat{s}u^\ell |^{n/2}d\mu \geq |Y_{[g]}|^{n/2},
$$
with equality iff $u$ is constant.
The result follows. 
\end{proof}

This leads to a very useful reinterpretation of $Y(M)$;
cf. \cite{and,bcg,okob,sch}. 

\begin{prop}\label{func} 
Let $M$ be a smooth compact $n$-manifold, $n \geq 3$. 
Then 
$$
\inf_g \int_M |s_g|^{n/2}d\mu_g = \left\{
\begin{array}{ll} 0&\mbox{if } Y(M)\geq 0\\
|Y(M)|^{n/2} &\mbox{if } Y(M)\leq 0.
\end{array}\right.
$$
Here the infimum on the left-hand side
is taken over all smooth Riemannian metrics
$g$ 
on $M$. 
\end{prop}
\begin{proof}
By Lemma \ref{yobo}, 
$$
\inf_g \int_M |s_g|^{n/2}d\mu_g = \inf_{[g]} |Y_{[g]}|^{n/2}.$$
If $Y(M) \leq 0$, we have $Y_{[g]} \leq 0$ for all $[g]$,
so the right-hand-side may be rewritten as 
$(-\sup Y_{[g]})^{n/2} = |Y(M)|^{n/2}$, and we are done. 
On the other hand, any
  smooth manifold of dimension $\geq 3$ admits metrics
of negative scalar curvature \cite{bes}, and 
$Y_{[g]}$  depends continuously
on $[g]$. Since    the space of metrics is 
connected,  $Y(M) > 0$ implies there is a conformal
class $[g]$ on $M$ with $Y_{[g]}=0$, and the infimum
therefore vanishes. 
\end{proof}

\section{Results from Seiberg-Witten Theory}

It is easy to see that 
some
complex surfaces of K\"ahler type
carry metrics with positive scalar curvature. 

\begin{prop} \label{rule}
Let $M$ be the underlying $4$-manifold
of a K\"ahler-type complex surface with 
Kodaira dimension $-\infty$. Then
$Y(M) > 0$. 
\end{prop}

Indeed, the Kodaira-Enriques classification \cite{bpv}
tells us that any such complex surface is either $\cp_2$ or else 
is obtained from a $\cp_1$-bundle over a complex
curve by blowing up points. Now $\cp_2$ carries the 
Fubini-Study metric, which has positive scalar curvature. 
Any $S^2$-bundle over another manifold also carries metrics
of positive scalar curvature; namely, consider any Riemannian
submersion metric with round, totally geodesic fibers, 
and then rescale the fibers to have very small diameter while
keeping
the metric on the base fixed. Finally, the blowing-up operation
amounts differentiably to taking connect sums with copies
of $\overline{\cp}_2$; all blow-ups of surfaces with 
positive-scalar-curvature metrics therefore 
themselves admit  
positive-scalar-curvature metrics because a result of 
Gromov and Lawson \cite{GL} tells us that the 
class of manifolds with positive-scalar-curvature metrics 
is closed under surgeries in codimension $\geq 3$.

On the other hand, it is a truly remarkable consequence of 
Seiberg-Witten theory that the converse of the
above statement is also true:
 
\begin{thm} \label{sw}
Let $M$ be the underlying 4-manifold of a 
complex surface of K\"ahler type with 
Kodaira dimension $\geq 0$. Then $Y(M) \leq 0$.
\end{thm} 

Remember, 
this is just a fancy way of saying that 
$M$ does not admit metrics of positive scalar curvature.
When $b^+(M) > 1$, this was first proved by
Witten \cite{witten}. In the $b^+=1$ case, this was
  proved for minimal surfaces in \cite{leb1}, and 
then in the   
non-minimal case by Friedman and Morgan \cite{FM}. 

The  proof proceeds by first
showing that for every metric
$g$ on $M$, there must be a solution of the Seiberg-Witten
equations 
\bea D^{\theta}\Phi &=&0 \\
 F_{\theta}^+&=&i \sigma(\Phi) \eea
and an unknown unitary connection $\theta$
on some line bundle $L\to M$  and 
 an unknown twisted spinor $\Phi\in \Gamma (S_+\otimes L^{1/2})$.
Here $D^\theta$ is the Dirac operator coupled to 
$\theta$, $ F_{\theta}^+$ is the projection
of the curvature of $\theta$ to the 
bundle $\Lambda^+$ of $g$-self-dual 2-forms, 
and $\sigma : S_+\otimes L^{1/2}\to \Lambda^+$
is the real-quadratic map induced by 
the natural 
isomorphism $\odot^2{\Bbb S}_+=\Lambda^+\otimes {\Bbb C}$. 
Moreover, $L$ can either be taken to be the anti-canonical line bundle
of $(M,J)$ or its pull-back via some 
diffeomorphism $M\to M$.  However, these equations imply
the Weitzenb\"ock formula
$$0=4\nabla^{\theta\ast}\nabla^\theta + s\Phi + |\Phi|^2\Phi,
$$
where $s$ is the scalar curvature of $g$. Taking the 
inner product with $\Phi$ and integrating, this tells
us
that 
$$ 0 = \int_M \left[ 
|2\nabla^\theta\Phi|^2+ s |\Phi |^2 + |\Phi |^4
\right]
d\mu ,$$
which is a contradiction if $s > 0$.

Pushing this argument further yields a stronger 
result \cite{leb3} 
for surfaces of general type:

\begin{thm}\label{gtype}
Let $M$ be the underlying 4-manifold of a 
complex surface with Kodaira dimension $2$. Then $Y(M) < 0$.
Moreover, if 
$X$ is the minimal model of $M$, then 
$Y(M)=Y(X)= -4\pi \sqrt{2c_1^2(X)}.$
\end{thm} 

(Recall \cite{bpv} 
that the minimal model $X$ is characterized by the 
fact that $M$ can be  blown down to $X$, but that $X$ cannot be 
blown down any further.  
The minimal model $X$ of a surface $M$ of general type
 is unique, and  $c_1^2  (X) > 0$.)

Indeed, applying the
 Schwarz inequality to the integrated Weitzenb\"ock formula  
  tells us that
$$(\int_M  s^2~d\mu )^{1/2} (\int_M
  |\Phi|^4 ~d\mu)^{1/2} \geq 
 \int_M  (-s)|\Phi|^2~d\mu  
\geq  \int_M  |\Phi|^4 ~d\mu .$$
We must therefore have 
$$\int_Ms^2d\mu \geq \int_M  |\Phi|^4 ~d\mu = 
8\int_M|F_{\theta}^+|^2d\mu$$
for any metric. However \cite{leb2}, 
$$ \int_M|F_{\theta}^+|^2d\mu \geq 4\pi^2([c_1(L)]^+)^2$$
where $[c_1(L)]^+$ is self-dual part of the 
$g$-harmonic deRham representative of $c_1(L)$.
Moreover, for each $g$ one can show \cite{leb3} that
 $L$ can be chosen so that $([c_1(L)]^+)^2\geq c_1^2(X)$. 
This shows that 
$$ \int_M s^2 d\mu \geq 32\pi^2 c_1^2(X)$$
for all metrics on $M$. But  by gluing
together K\"ahler-Einstein orbifolds, gravitional instantons, and
standard metrics  on the blow-up of ${\Bbb C}^2$ at the
origin, one can also 
construct \cite{leb3} 
 sequences of metrics  
which show that this bound is actually sharp:
$$\inf_g \int_M s^2 d\mu = 32\pi^2 c_1^2(X) >0.$$
Thus $Y(M) = -4\pi \sqrt{2c_1^2(X)} < 0$ by
Proposition \ref{func}.

\section{Rational Elliptic Surfaces Revisited}

A complex surface $(M^4,J)$ is said to be {\em elliptic} if it admits
a holomorphic map $M\to \Sigma$ to a complex curve such that
the regular fibers are elliptic curves $T^2$. While this notion
is of greatest importance in the study of surfaces of
Kodaira dimension $0$ and $1$, there are also elliptic 
surfaces of Kodaira dimension $-\infty$. One family of 
examples of the latter type are the rational elliptic surfaces
treated this section. 

The smooth oriented 4-manifold $\cp_2\# 9 \overline{\cp}_2$
may be concretely realized as the blow-up of $\cp_2$
at $9$ points. One  choice of these $9$ points is
to take them to be the points of intersection of a generic pair of
elliptic (i.e. cubic) curves in the projective plane. 
There is then exactly  a $\cp_1$'s worth of cubics passing through the 
$9$ points, corresponding to linear combinations of the 
 defining equations of the original pair of curves. 
When  $\cp_2$ is blown up at these  $9$   points,
this pencil of cubics  gives rise to a holomorphic projection
 $\cp_2\# 9 \overline{\cp}_2\to \cp_1$, the generic fiber
of which is a smooth elliptic curve $E\approx T^2$. 
The term {\em rational elliptic surface} is  thus 
used to describe    
$\cp_2\# 9 \overline{\cp}_2$ equipped
 with  the deformation class of this projection to 
$\cp_1$.

There is a  useful way of constructing special 
rational elliptic surfaces by analogy to 
the Kummer construction of K3's. 
Let $\Lambda \cong {\Bbb Z}\oplus {\Bbb Z}$
be any lattice in $\Bbb C$, and let
$E={\Bbb C}/\Lambda$ be the corresponding elliptic curve. 
Recall that
each such curve can be realized as a branched cover over
$\cp_1$ because reflection $z\mapsto -z$
through the origin in $\Bbb C$ induces in an involution
of $E$ with $4$ fixed-points\footnote{corresponding the
half-lattice $\frac{1}{2}\Lambda$}; the orbifold quotient 
$E/{\Bbb Z}_2$ can then be identified with 
a Riemann surface of genus $0$, and this 
realizes $E$ as a branched cover   
$\wp : E\to \bcp_1$. The involution of $E$ that
occurs in this argument is called the {\em Weierstrass 
involution} of $E$, and will be denoted here by
$\Phi : E\to E$. Notice that $\Phi$ is an isometry
of the obvious flat metric on $E$.

Now
let $\Psi : \cp_1 \to \cp_1$ be the involution
$[z_0:z_1]\mapsto [-z_0: z_1]$, which has two fixed
points and amounts to a $180^\circ$
rotation of $S^2$ around an axis. The product
involution 
$$\Phi \times \Psi : E\times \cp_1 \to E\times \cp_1$$
then 
has 8 fixed points, and 
the quotient $E\times \cp_1/{\Bbb Z}_2$ is
a complex orbifold with 8 singularities
modeled on ${\Bbb C}^2/\pm 1$. Let $\tilde{M}$ be
the blow-up of 
of $E\times \cp_1$ at the 8 fixed points of 
$\Phi \times \Psi$, and let $M = \tilde{M}/{\Bbb Z}_2$.
The non-singular complex surface $M$ is thus
obtained  from $E\times \cp_1/{\Bbb Z}_2$ 
by replacing each of the $8$ singular points with a
$(-2)$-curve.

Now consider the  holomorphic map
$\varpi : M \to \cp_1$ 
induced by following the
second-factor projection $E\times \cp_1\to \cp_1$
with $[z_0: z_1] \mapsto [z_0^2 : z_1^2]$. All
but two fibers of $\varpi$ are copies of $E$;
the only two exceptions are the fibers over the
north and south poles $[0:1]$ and $[1:0]$. Each of these 
two exceptional fibers  consists of five $(-2)$-curves, 
linked together according to the Dynkin diagram 
$\tilde{D}_4$:

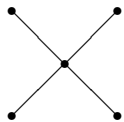
\begin{center}
\begin{picture}(240,60)(0,3)
\put(100,30){\circle*{3}}
\put(80,50){\circle*{3}}
\put(80,10){\circle*{3}}
\put(120,50){\circle*{3}}
\put(120,10){\circle*{3}}
 \put(80,10){\line(1,1){40}}
\put(80,50){\line(1,-1){40}}
\end{picture}
\end{center}
These  two exceptional
fibers are of thus of Kodaira type $I_0^*(\tilde{D}_4)$.

Now  
$\varpi : M \to \bcp_1$ is deformation
equivalent to the above-described standard model of the
rational elliptic surface.  
In fact, $M$ can be blown down to $\cp_1\times \cp_1$
in such a way that the fibers of $\varpi$ are sent to
the elliptic curves of the pencil
$$t_0x_0(x_0-x_1)y_0^2+ t_1x_1(x_0-ax_1)y_1^2=0,$$
where the number $a$ is given in terms of the lattice
$\Lambda$ by the Weierstrass  $\jmath$-function.  
Notice that for each $[t]\in \cp_1$ this locus is
 an element
of the anti-canonical linear system $|{\cal O}(2,2)|$ on
$\cp_1\times\cp_1$, and that each  
elliptic curve of the pencil passes through the
points    $([0:1],[0:1],)$,
$([1:1],[0:1])$,
$([a:1],[1:0])$, and 
$([1:0],[1:0])$, and all are 
 tangent to the 
  second-factor $\cp_1$ at these points.

To see this explicitly, let us first
consider $E\times \cp_1$ as a ruled surface
over $E$, and blow it up at the $8$ 
relevant points to obtain $\tilde{M}$:
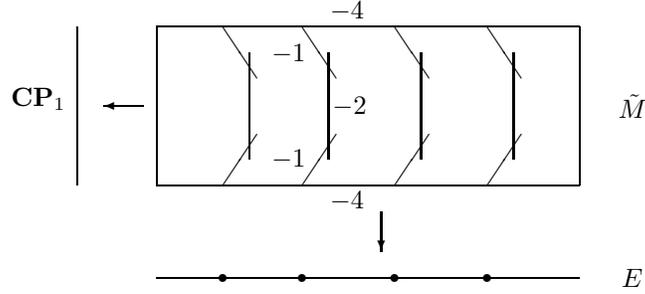
\begin{center}
\begin{picture}(240,110)(0,3)
\put(50,110){\line(1,0){160}}
\put(50,110){\line(0,-1){60}}

\put(20,110){\line(0,-1){60}}
\put(45,80){\vector(-1,0){15}}
\put(-5,80){$\cp_1$}

\put(75,110){\line(2,-3){13}}
\put(85,100){\line(0,-1){40}}
\put(75,50){\line(2,3){13}}

\put(105,110){\line(2,-3){13}}
\put(115,100){\line(0,-1){40}}
\put(105,50){\line(2,3){13}}
\put(100,100){\makebox(0,0){$-1$}}
\put(100,60){\makebox(0,0){$-1$}}
\put(123,80){\makebox(0,0){$-2$}}

\put(140,110){\line(2,-3){13}}
\put(150,100){\line(0,-1){40}}
\put(140,50){\line(2,3){13}}

\put(175,110){\line(2,-3){13}}
\put(185,100){\line(0,-1){40}}
\put(175,50){\line(2,3){13}}

\put(122,44){\makebox(0,0){$-4$}}
\put(122,116){\makebox(0,0){$-4$}}
\put(230,80){\makebox(0,0){$\tilde{M}$}}
\put(230,15){\makebox(0,0){$E$}}

\put(50,50){\line(1,0){160}}
\put(210,110){\line(0,-1){60}}
\put(135,40){\vector(0,-1){15}}
\put(50,15){\line(1,0){160}}

\put(75,15){\circle*{3}}
\put(105,15){\circle*{3}}
\put(140,15){\circle*{3}}
\put(175,15){\circle*{3}}
\end{picture}
\end{center}
Dividing by $\Phi\times\Psi$ then gives us a
non-minimal rational ruled surface:
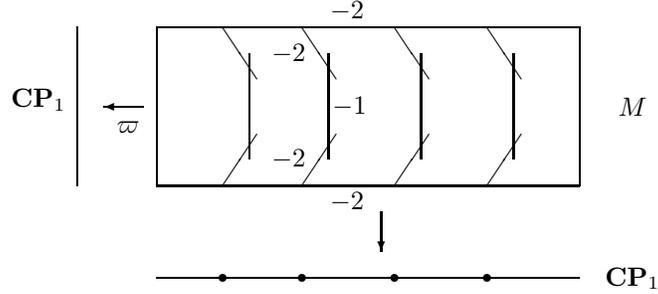
\begin{center}
\begin{picture}(240,110)(0,3)
\put(50,110){\line(1,0){160}}
\put(50,110){\line(0,-1){60}}

\put(20,110){\line(0,-1){60}}
\put(45,80){\vector(-1,0){15}}
\put(35,70){$\varpi$}
\put(-5,80){$\cp_1$}

\put(75,110){\line(2,-3){13}}
\put(85,100){\line(0,-1){40}}
\put(75,50){\line(2,3){13}}

\put(105,110){\line(2,-3){13}}
\put(115,100){\line(0,-1){40}}
\put(105,50){\line(2,3){13}}
\put(100,100){\makebox(0,0){$-2$}}
\put(100,60){\makebox(0,0){$-2$}}
\put(123,80){\makebox(0,0){$-1$}}

\put(140,110){\line(2,-3){13}}
\put(150,100){\line(0,-1){40}}
\put(140,50){\line(2,3){13}}

\put(175,110){\line(2,-3){13}}
\put(185,100){\line(0,-1){40}}
\put(175,50){\line(2,3){13}}

\put(122,44){\makebox(0,0){$-2$}}
\put(122,116){\makebox(0,0){$-2$}}
\put(230,80){\makebox(0,0){${M}$}}
\put(230,15){\makebox(0,0){$\cp_1$}}

\put(50,50){\line(1,0){160}}
\put(210,110){\line(0,-1){60}}
\put(135,40){\vector(0,-1){15}}
\put(50,15){\line(1,0){160}}

\put(75,15){\circle*{3}}
\put(105,15){\circle*{3}}
\put(140,15){\circle*{3}}
\put(175,15){\circle*{3}}
\end{picture}
\end{center}
(Note that 
$\varpi$, which is depicted horizontally in the 
above picture, is {\em not}
the ruled-surface projection, which is 
here 
depicted vertically.)  
Contracting  $8$ exceptional  curves in a judicious order,

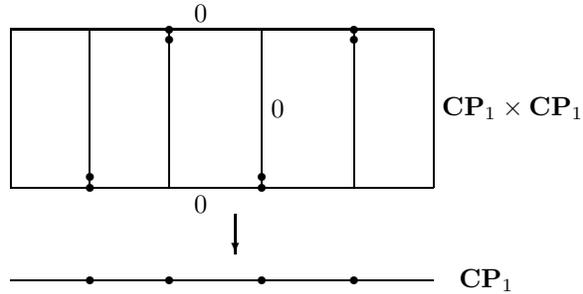
\begin{center}
\begin{picture}(240,110)(0,3)
\put(50,110){\line(1,0){160}}
\put(50,110){\line(0,-1){60}}

\put(80,110){\line(0,-1){60}}
\put(110,110){\line(0,-1){60}}
\put(145,110){\line(0,-1){60}}
\put(180,110){\line(0,-1){60}}

\put(80,50){\circle*{3}}
\put(145,50){\circle*{3}}
 
\put(80,54){\circle*{3}}
\put(145,54){\circle*{3}}

\put(110,110){\circle*{3}}
\put(180,110){\circle*{3}}
 
\put(110,106){\circle*{3}}
\put(180,106){\circle*{3}}

\put(122,44){\makebox(0,0){$0$}}
\put(122,116){\makebox(0,0){$0$}}
\put(240,80){\makebox(0,0){$\cp_1 \times \cp_1$}}
\put(230,15){\makebox(0,0){$\cp_1$}}

\put(50,50){\line(1,0){160}}
\put(210,110){\line(0,-1){60}}
\put(135,40){\vector(0,-1){15}}
\put(50,15){\line(1,0){160}}

\put(151,80){\makebox(0,0){$0$}}

\put(80,15){\circle*{3}}
\put(110,15){\circle*{3}}
\put(145,15){\circle*{3}}
\put(180,15){\circle*{3}}
\end{picture}
\end{center}
we get a rational ruled surface which can be 
recognized as  
$\bcp_1\times\bcp_1$ by inspecting self-intersections
of holomorphic curves. Our   pencil now becomes 
the sub-family of the linear system $|{\cal O}(2,2)|$
consisting of curves which are tangent to the second factor at 
$4$   points in the claimed special position. 
It follows that $\varpi : M \to \cp_1$ is 
deformation equivalent to the 
standard model of the rational elliptic surface,
since we may first deform our pencil into a generic
sub-pencil of $|{\cal O}(2,2)|$, and then 
map $\cp_1\times\cp_1$ birationally
to $\cp_2$ in the usual way, using one 
of the
base points of the pencil as our center. 

This leads to the following:

\begin{prop} \label{eh}
Let $V$ be the (flat) 4-orbifold obtained from 
${\Bbb R}\times T^3$ by dividing by the
involution induced by  
$-1:  {\Bbb R}^4\to {\Bbb R}^4$,   
let $B=({\Bbb R}\times S^1)/{\Bbb Z}_2$ denote the analogous 
2-orbifold, and let $\pi : V\to B$ be the 
map induced by projection ${\Bbb R}\times T^3\to {\Bbb R}\times S^1$
to the first two coordinates. Then 
a smooth model of the  rational elliptic
surface $\cp_2\# 9\overline{\cp}_2 \to \cp_1$
may be obtained from $\pi : V\to B$ by 
replacing each of the $8$ singular point of $V$ with a
2-sphere  of self-intersection $-2$, forgetting the 
  orbifold structure on $B$, and adding 
a smooth fiber at infinity in   the
 obvious manner. 
\end{prop}
 \begin{proof}
We have already seen that the rational elliptic
surface can be viewed as a desingularization 
of $(S^2\times T^2)/{\Bbb Z}_2$. Now let
$x\in S^2$ be some point which is not fixed by 
the involution, and let $x^\prime\neq x$ 
be its image under the involution. 
Identify $S^2-\{x,x^\prime\}$ with 
the cylinder ${\Bbb R}\times S^1$ in such a 
manner that the involution becomes 
simultaneous reflection in both factors. The result follows.
\end{proof}

The orbifold $V$ comes equipped with a family of
flat orbifold metrics. Replacing the 
singularities with gravitational instantons yields
the main result of this section:

\begin{prop} Let $\tilde{V}$ denote the complement of 
a generic fiber in the rational elliptic
surface $\cp_2\# 9 \overline{\cp}_2\to \cp_1$,  
let $\pi: \tilde{V}\to {\Bbb C}$ be the induced elliptic
fibration, normalized so that all critical 
values are contained in open unit disk    $\Delta\subset {\Bbb C}$.
 Let $f$ be any given flat  metric on 
the 2-torus $T^2$. 
 Then there is a family $g_t$, $t\in [1,\infty )$ 
of Riemannian metrics
 on $\tilde{V}$ 
such that 
\begin{itemize}
\item the Ricci curvature of $g_t$ converges uniformly to $0$ as 
$t\to\infty$; 
\item $(\pi^{-1} ({\Bbb C} - \Delta),g_t)$
is isometric to $([0,\infty )\times S^1, dx^2+d\theta^2)
\times (T^2, f/t)$; and 
\item $\lim_{j\to \infty} \mbox{ Vol } (\pi^{-1}(\Delta), g_t ) = 0$. 
\end{itemize}
\end{prop}
\begin{proof}
The basic idea is to replace the singularities of
of $V$ with  Eguchi-Hanson metrics.
The {\em Eguchi-Hanson metric} \cite{eh,leb0} is a complete Ricci-flat 
 metric on the ${\cal O}(-2)$ line bundle over $\cp_1$,
or in other words on the blow-up of ${\Bbb C}^2/(\pm 1)$ 
at the origin.
Essentially by introducing polar coordinates on 
 $({\Bbb R}^4-0)/(\pm 1)$, we may think of this as 
the completion of the metric on 
 $(1,\infty)\times SO(3)$ given by   
$$g_{EH}=\frac{dr^2}{1-\frac{1}{r^4}}+ r^2\left( \sigma_1^2+\sigma_2^2+ 
(1-\frac{1}{r^4})\sigma_3^2\right)$$
where the left-invariant co-frame $\{ \sigma_j\}$ is
orthonormal for the curvature $+1$ bi-invariant metric. 
By rescaling and homothety, we can turn this into
a 1-parameter family 
$$\frac{dr^2}{1-\frac{A}{r^4}}+ r^2\left( \sigma_1^2+\sigma_2^2+ 
(1-\frac{A}{r^4})\sigma_3^2\right), ~~ r > \sqrt[4]{A}$$
of Ricci-flat metrics which bubble off to
the flat orbifold metric on ${\Bbb R}^4/(\pm 1)$
as $A\to 0$.

Let us now modify the Eguchi-Hanson metric. Let 
$\phi (t)$ be a smooth monotonely non-increasing function 
which is identically $1$ on
$[0,1]$ and  identically $0$ on $[2,\infty ]$.
For small values of the parameter $A$, 
the modified metric 
$$\frac{dr^2}{1-\phi(r)\frac{A}{r^4}}+ r^2\left( \sigma_1^2+\sigma_2^2+ 
(1-\phi(r)\frac{A}{r^4})\sigma_3^2\right), ~~r > \sqrt[4]{A}$$
then agrees with the Euclidean metric for
$r > 2$, is Ricci-flat when $r < 1$, and has 
Ricci curvature $< CA$ for some constant $C$ independent
of the small parameter $A$. Making 
another homothety and rescaling, the metric 
$$\frac{dr^2}{1-\phi(r/\varepsilon)\frac{A\varepsilon^4}{r^4}}+ 
r^2\left( \sigma_1^2+\sigma_2^2+ 
(1-\phi(r/\varepsilon)\frac{A\varepsilon^4}{r^4})\sigma_3^2\right), 
~~r > \varepsilon \sqrt[4]{A}$$
is Euclidean for $r > 2\varepsilon$, Ricci-flat for $r < \varepsilon$,
and has Ricci curvature $< CA/\varepsilon^2$. Thus the metrics
$$g_\varepsilon=\frac{dr^2}{1-\phi(r/\varepsilon)
\frac{\varepsilon^8}{r^4}}+ 
r^2\left( \sigma_1^2+\sigma_2^2+ 
(1-\phi(r/\varepsilon)\frac{\varepsilon^8}{r^4})\sigma_3^2\right),
~~ r> \varepsilon^2$$
 have    Ricci curvature uniformly $O(\varepsilon^2)$
as $\varepsilon \searrow 0$.
Notice that the volume form of $g_\varepsilon$ coincides with the
Euclidean volume form 
$r^3dr\wedge\sigma_1\wedge\sigma_2\wedge\sigma_3$,
but that the appropriate domain of integration is 
the   $r > \varepsilon^2$
rather than $r > 0$. Thus,
  smoothing a flat orbifold   singularity modeled 
on ${\Bbb R}^4/(\pm 1)$ by gluing in the
modified Eguchi-Hanson  metric $g_\varepsilon$
 reduces the  total volume, namely by $\pi^2\varepsilon^8/2$.

Now let $f$ be any fixed flat metric
on the 2-torus, and let $\imath$ denote the injectivity 
radius of $(T^2,f)$, and
set $\hat{\imath}=\min (\imath , \pi)$.
 Endow ${\Bbb R}\times T^3$ with the sequence
of  flat metrics
$\hat{g_t}=dx^2+d\theta^2+ f/t$, and push these
down as flat orbifold metrics on $V=({\Bbb R}\times T^3)/{\Bbb Z}_2$.
For each $t$, let $\varepsilon_t = \hat{\imath}/4\sqrt{t}$, 
so that the balls of radius $2\varepsilon_t$ centered at the
orbifold singularities of $V$ are pairwise disjoint. 
On the blow-up $\tilde{V}$ of $V$, we can then define a 
metric $g_t$ as the $\hat{g_t}$ on the complement of
these balls, and the modified Eguchi-Hanson metrics
$g_{\varepsilon_t}$ on the blown-up interior of these balls . 
Since $\varepsilon_t\to 0$, the Ricci curvature of $g_t$
tends uniformly to zero. Moreover,  the volume of
the region $|x| \leq a$ in $(X,g_t)$ is less than 
$2\pi a \alpha/t \to 0$, where $\alpha$ is the area of 
$(T^2,f)$, so the volume of the pre-image of any 
compact set in $\Bbb C$ tends to $0$ as $t\to\infty$. 
\end{proof}

\section{Collapsing  Elliptic Surfaces}

\begin{defn}
Let $M$ be a smooth manifold. We will say that $M$
{\em collapses with bounded scalar curvature}
if there is a sequence $g_j$ of smooth metrics 
on $M$ for which the absolute value  
of the scalar curvatures is uniformly bounded, 
but with total volume tending to zero:
$$|s_{g_j}| < {\cal B}, ~~\lim_{j\to \infty}\mbox{Vol}(M,g_j) =0.$$
Similarly, we say that $M$ collapses with 
bounded Ricci (respectively,   sectional) 
curvature  if there is  a sequence $g_j$ of smooth metrics
on $M$ with 
  uniformly  bounded Ricci (respectively, sectional)
curvature and total volume tending to zero. 
\end{defn}

This collection of definitions is loosely inspired by the 
work of Cheeger and Gromov \cite{cg}, who studied 
sequences of metrics $g_j$ with bounded sectional curvature
and {\em injectivity radius} tending to zero at all points. 
A standard comparison argument for the volume of 
small balls shows that 
 a sequence of metrics with sectional   
curvature bounded above can have volume tending with to zero 
only if  the point-wise injectivity radii 
uniformly tend to zero, too. 
Note, however, that a   sequence
of metrics $g_j$ with bounded sectional curvature 
and injectivity radius tending to zero may have volume
bounded away from zero --- as  even happens for 
flat  metrics on the 2-torus. 
Nonetheless, a remarkable result of Cheeger and Gromov \cite{cg2}
asserts 
that if a compact manifold admits   
 a sequence of metrics $g_j$
with bounded sectional curvature and injectivity radius
uniformly tending to zero,  
 there must be other sequences
$\tilde{g}_j$ with bounded sectional curvature for which the 
volume tends to zero as well.  Thus the
present definition of collapse with bounded sectional
curvature  
is in accord with the definitions used by others.

Collapse with bounded scalar curvature 
is directly relevant to the computation
of Yamabe invariants  by virtue of the
following result:

\begin{prop} \label{sam} 
Let $M$ be a smooth compact $n$-manifold, $n\geq 3$. Then
the following are equivalent:
\begin{description}
\item{(i)} $M$ collapses with bounded scalar curvature;
\item{(ii)} $\inf_{g}\int_M |s_g|^{n/2}d\mu_g =0$; 
\item{(iii)} $Y(M)\geq 0$.
\end{description}
\end{prop}
\begin{proof}
The implication $(i)\Rightarrow (ii)$ is trivial, since
$\int_M |s_{g_j}|^{n/2}d\mu_{g_j}\to 0$ for any sequence of
metrics with bounded scalar curvature and volume tending to 
zero. The implication $(ii)\Rightarrow (iii)$ follows
from Proposition \ref{func}. Finally, 
the  implication $(iii)\Rightarrow (i)$ 
follows from the fact \cite{aubin} that any $n$-manifold, $n \geq 3$,
admits of metrics of constant
negative scalar curvature; if $Y(M)\geq 0$, 
 there is a thus a sequence of unit-volume metrics with
constant negative scalar curvature tending to $0$ from below,
and rescaling these metrics so they have $s\equiv -1$ then
makes the volumes of the rescaled metrics tend to zero. 
\end{proof}

If $M$ is the underlying 4-manifold of a
complex  elliptic surface, these various notions of 
collapse can be completely understood.

\begin{thm}\label{scol}
Let $M$ be the underlying 4-manifold of a
complex elliptic surface. Then $M$ collapses with bounded
scalar curvature.
\end{thm}

\begin{thm}\label{rcol}
Let $M$ be the underlying 4-manifold of a
complex elliptic surface. Then $M$ collapses with bounded
Ricci curvature $\Leftrightarrow$ 
$M$ is minimal. 
\end{thm}

\begin{thm}\label{kcol}
Let $M$ be the underlying 4-manifold of a
complex elliptic surface. Then $M$ collapses with bounded
sectional curvature $\Leftrightarrow$ 
$\chi(M)=0$.
\end{thm}

In order to prove these results, we will 
need special diffeomorphic models \cite{FM,mats} of elliptic
surfaces. These are built up in stages, starting with 
surfaces we shall call {\em twisted products}. 
Let $E$ be an elliptic curve, equipped with a compatible flat
metric, let $o\in E$ be a chosen base-point, 
and let $G$ be the finite group of 
orientation-preserving isometries $E\to E$ which 
fix $o$. (Thus $G$ is a cyclic group of order
$2$, $4$, or $6$.) A twisted product surface
is simply the total space $B$ of a fiber-bundle $B\to \Sigma$ over a 
compact complex curve, with fiber $E$ and
structure group $G$.

 Now let $B\to \Sigma$ be such a twisted product, 
and let a finite collection of points $x_j\in \Sigma$ 
be specified, and assign each of these points 
 some integer multiplicity
$m_j > 1$.  Give
$\Sigma$ an obifold structure by introducinging orbifold
local coordinates  of the form $\hat{z}={z}^{1/m_j}$
near $x_j$, where
$z$ is a local complex coordinate on $\Sigma$   with $z(x_j)=0$. 
A process called  the {\em logarithmic transform}
now allows us to  delete the fiber $E_j\cong E$ over each $x_j$, 
and replace
it with $E_j/{\Bbb Z}_{m_j}$, where the ${\Bbb Z}_{m_j}$-action
is generated by any specified translation 
$\tau_j : E_j\to E_j$ of order $m_j$. Namely, let 
$\tau_j(t)$ be a 1-complex-parameter group of translations
of $E_{x_j}$
with $\tau_j (1)=\tau_j$, let $\Delta\subset {\Bbb C}$ be a  disk centered 
at $0$ and contained in the range of the orbifold coordinate
$\hat{z}$, and glue $(E_j\times \Delta)/\langle (\tau_j, e^{2\pi i/m_j})\rangle$
to $B-E_{x_j}$ by 
$[(y,\tilde{z})]\leftrightarrow (\tau_j (-m_j\log (\hat{z})/2\pi i )(y),
z:=\hat{z}^{m_j})$. 

With these constructions in mind, we can now
prove

\begin{lem} Let $M$ be the underlying 4-manifold of any 
 complex elliptic surface  with Euler characteristic
$0$. Then there is a  compact $2$-orbifold   $\Sigma$
and a
smooth submersive map $\pi : M \to \Sigma$, 
such that for every orbifold metric $h$ on $\Sigma$
there is 
a Riemannian metric $g$ on $M$ for which 
$\pi$ becomes a Riemannian submersion
with  flat, totally geodesic 
fibers. Moreover, $\pi$ can be chosen so
that its fibers over some open disk in $\Sigma$ are also
regular fibers of the original elliptic fibration of $M$. 
\end{lem}
\begin{proof} Any complex elliptic surface 
with Euler characteristic zero is 
$C^\infty$-deformation equivalent to some logarithmic  transform of 
a twisted-product surface  \cite[Proposition 7.2]{FM}. 
But this model is an orbi-bundle over the $2$-orbifold $\Sigma$
with structure group equal to the isometry group of
a flat $2$-torus $E$. Starting from any 
orbifold metric $h$ on $\Sigma$, now construct a metric
$g$ on $M$ by gluing together local product metrics
with a  partition of unity subordinate to a trivializing 
cover of $\Sigma$. Each fiber of $\pi: M\to \Sigma$
is then a flat quotient of $E$, and each fiber has vanishing
second fundamental form because  the orthogonal space of
the fibers is the horizontal subspace of a connection
with structure group equal to the isometry group of the fiber. 
\end{proof}

The ideas of Cheeger and Gromov \cite{cg}
now immediately lead to a 

\bigskip

\noindent {\bf Proof of Theorem \ref{kcol}.}
Consider the metrics
$$g_t = \frac{1}{t}g+(1-\frac{1}{t})\pi^*h$$
on $M$ as $t\to\infty$. We have
$$\mbox{Vol}(M,g_t)= \frac{1}{t}\mbox{Vol}(M,g)\to 0,$$
so it suffices to show that the components of 
the curvature tensor remain bounded in an orthonormal
frame as $t\to\infty$.  But this follows from O'Neill's
 Riemannian submersion formulas \cite{bes}. Indeed, the vertical
sectional curvature is $0$ for all $t$, and the
horizontal sectional curvature is given by
$$K(H)=K(\Sigma) - \frac{3}{4}g_t(v,v)\to K(\Sigma),$$
while the sectional curvature of the 2-plane $P$
spanned by a vertical vector and
a horizontal vector is
$$K(P)= \frac{1}{4}g_t(v^{\|},v^{\|}) \to 0.$$
Here $K(\Sigma)$ is 
the Gauss curvature of $(\Sigma,h)$, and 
$v$ is the ($t$-independent) vertical component of $[w_1,w_2]$, 
where $w_1$ and $w_2$ are horizontal vectors which 
project to an oriented orthonormal frame on $\Sigma$,
whereas $v^{\|}$ is the 
orthogonal projection of $v$ into $P$.   
This shows that the point-wise norm of the curvature tensor,
 and hence all sectional curvatures,
 remain uniformly bounded
as $t\to \infty$. Any elliptic surface $M$ 
with $\chi (M)=0$ therefore collapses 
with bounded sectional curvature. 

Conversely, if $M$  collapses with bounded sectional 
curvature, 
the  Gauss-Bonnet formula tells one that
$\chi (M) =0$.
\qed

The next operation we will need in order to 
create smooth models of elliptic surfaces is the 
{\em fiber sum}. The fiber sum of two elliptic
surfaces is the smooth manifold obtained by 
removing a tubular neighborhood of a regular fiber in 
each and then then identifying the boundaries 
in a manner compatible with the orientations and
given local trivializations. 
Any minimal complex elliptic surface  $M$ is
 diffeomorphic \cite[Corollary 2.17]{FM} to a fiber sum of  
an elliptic surface  $\check{M}$ 
with Euler characteristic $0$
and copies of the rational elliptic surface. 

\bigskip

\noindent {\bf Proof of Theorem \ref{rcol}.}
  Let us thus begin with some elliptic surface $\check{M}$
with $\chi (\check{M})=0$, and let $\pi : \check{M}\to \Sigma$
be a smooth submersion of the type used in the previous proof. 
Let $x_1, \ldots , x_k\in \Sigma$ be the points we will use as the
centers of the fiber sum construction. Let $U_1, \ldots , U_k$
be disjoint disks around these points, chosen so that
$M\to \Sigma$ may be trivialized over their closures, 
and identify
each with the   disk of radius $2$ about   $0\in {\Bbb C}$. Choose $h$ so
that its restriction to the annulus $|z|\in [ 1, 2]$ 
is isometric to the cylinder $S^1\times [1,2]$, and now
construct the family of metrics $g_t$ used above, with 
the stipulation $g$ is actually taken to be a product metric on
$\pi^{-1}(U_1), \ldots , \pi^{-1} (U_k)$.
Letting $\tilde{U}_1\subset U_1, \ldots , \tilde{U}_k\subset U_k$ 
each correspond to
the unit sub-disk $|z| < 1$,  
the Riemannian 
manifold-with-boundary $(M-\cup \pi^{-1}(\tilde{U}_j),g_t)$ has ends
isometric to $([0,1]\times S^1, dt^2+ d\theta^2) \times (E, f/t)$,
where $f$ is a fixed flat metric on an elliptic curve $E$. 

However, these cylindrical ends
 precisely match those of the metrics on $\tilde{V}$ 
constructed in Proposition \ref{eh}, and we can therefore
glue these   objects together to obtain metrics with 
bounded Ricci curvature and volume tending to
zero on the fiber sum $M$ of $\check{M}$ and $k$
rational elliptic surfaces.  Since a diffeomorphic
model of any minimal elliptic surface can be 
constructed in this way, it follows that  any minimal 
elliptic surface collapses
with bounded Ricci curvature, as claimed.

Conversely, suppose
 an oriented 4-manifold $M$ admits
 a sequence of metrics $g_j$ with bounded Ricci
curvature $r$ and volume tending to zero. We then 
have
$$\lim_{j\to\infty} \int_M \left[ \frac{s_{g_j}^2}{24}-
\frac{|\stackrel{\circ}{r}|^2_{g_j}}{2}\right]d\mu_{g_j}
=0,$$
where $\stackrel{\circ}{r}$ denotes the 
trace-free part of the Ricci curvature. 
But 
$$(2\chi + 3\tau)(M)= \frac{1}{4\pi^2}\int_M \left[ |W_+|^2+
\frac{s_{g}^2}{24}-
\frac{|\stackrel{\circ}{r}|^2_{g}}{2}\right]d\mu_{g}
\geq \frac{1}{4\pi^2}\int_M \left[ \frac{s_{g}^2}{24}-
\frac{|\stackrel{\circ}{r}|^2_{g}}{2}\right]d\mu_{g}$$
for any metric on $M$. It therefore follows that 
$(2\chi + 3\tau )(M) \geq 0$. 
On the other hand, $2\chi + 3\tau = c_1^2 \leq 0$ if
$M$ is an elliptic surface, and equality occurs iff
$M$ is minimal \cite{bpv}. Hence an elliptic
surface $M$ must be
minimal if it collapses with bounded Ricci curvature. 
\qed

Finally, any elliptic surface is obtained from
a minimal elliptic surface by {\em blowing up} 
an appropriate number of times. 
Recall that the blow-up of ${\Bbb C}^2$ at the
origin is by definition the line bundle ${\cal O}(-1)$
over $\cp_1$.  A point in an arbitrary complex surface can
similarly be  
obtained by replacing a small ball with a 
tubular neighborhood of
the zero section in this bundle. This procedure  
is then iterated as needed to produce the desired 
surface. 
\bigskip

\noindent {\bf Proof of Theorem \ref{scol}.}
The idea is to graft Burns metrics  onto the previous 
examples. The {\em Burns metric} \cite{leb0} is a scalar-flat
K\"ahler metric on 
the blow-up of ${\Bbb C}^2$ at the origin.
By introducing  polar coordinates, this can be viewed
as the completion of the metric on $(1,\infty)\times S^3$
given by 
$$g_{\mbox{\tiny 
Burns}}= \frac{dr^2}{1-\frac{1}{r^2}} + r^2\left(\sigma_1^2+
\sigma_2^2+ (1-\frac{1}{r^2})\sigma_3^2\right) ,$$
where $\{\sigma_j\}$ is an $SU(2)$-invariant orthonormal
frame on $S^3$. By repeating the sequence of homotheties, 
rescalings, and cut-offs used in the proof of
Proposition \ref{eh}, the metric 
$$g_\epsilon = 
\frac{dr^2}{1-\phi(\frac{r}{\epsilon})\frac{\epsilon^6}{r^2}} 
+ r^2\left(\sigma_1^2+
\sigma_2^2+ (1-\phi(\frac{r}{\epsilon})\frac{\epsilon^6}{r^2}
)\sigma_3^2\right),~~
r > \epsilon^3  
$$
is flat for $r > 2\epsilon$, scalar-flat for $r< \epsilon$,
and has scalar curvature    uniformly $O(\epsilon^2)$
as $\epsilon \searrow 0$. Moreover, blowing up a
flat region by replacing a small ball by a
copy of $g_\epsilon$ reduces its total volume, namely
by $\pi^2\epsilon^{12}/2$.

Any elliptic surface is an iterated blow-up of a minimal
one. Up to diffeomorphism, we can therefore 
model any such surface by blowing up one of 
our previous models at distinct points. But the
constructed metrics on the previous models 
can   be taken to be flat in the neighborhood
of certain fibers  
provided the orbifold metric 
$h$ on $\Sigma$ is also taken to be flat on certain regions.
So let $\check{g}_j$ be a sequence of metrics with 
bounded Ricci curvature   and Vol $\to 0$,
each containing a Euclidean ball of radius $\varrho_j$,
on a minimal elliptic surface $\check{M}$.
On any
blow-up $M=\check{M}\# \ell \overline{\cp}_2$,
 we can then 
produce a sequence of metrics $g_j$, 
with bounded scalar curvature and Vol $\to 0$
by replacing $\ell$ disjoint, Euclidean $(2{\epsilon_j})$-balls with 
copies of $g_{\epsilon_j}$, where we may, for example,
take $\epsilon_j=\varrho_j/2\ell$. 
\qed

By Proposition \ref{sam}, this immediately implies 

\begin{cor}
Any complex elliptic surface $M$ has $Y(M) \geq 0$.
\end{cor}

The reader may prefer to deduce this corollary
directly from  Theorem 
\ref{rcol} and the 
fact that $Y(\overline{\cp}_2) > 0$ by 
citing a general result of \cite{okob}.
At heart, however, this is much the 
same proof as is given above. 

Combining this corollary with Theorem \ref{sw} now yields 

\begin{thm} \label{y0}
Let $M$ be a K\"ahler-type complex    surface   
with $\kod (M)=0$ or  $1$. Then    $Y(M)=0$. 
\end{thm}
\begin{proof}
Any surface of Kodaira dimension $0$ or $1$ is 
deformation-equivalent to an elliptic surface \cite{bpv,FM},
and so has $Y(M)\geq 0$ by the previous result. 
But if $b_1(M)=0\bmod 2$, $M$ is of K\"ahler type,
and   Seiberg-Witten theory tells us that 
there is no metric of positive scalar curvature 
on $M$, so that $Y(M)\leq 0$. Hence $Y(M)=0$, as claimed.
\end{proof}

Combining 
this  with Proposition
\ref{rule} and Theorem \ref{gtype}  
now proves Theorem \ref{main}.

\section{Concluding Remarks}

Theorem \ref{main} does not allow one to distinguish 
between Kodaira dimensions $0$ and $1$ by 
comparing Yamabe 
invariants. Such a  distinction arises 
immediately, however, if one asks whether
$$Y(M)=\sup_{[g]} Y_{[g]}$$
is actually realized by the Yamabe constant
 of some conformal class.

\begin{prop}
Let $M$ be the underlying 4-manifold of
a complex algebraic surface of Kodaira dimension
$\geq 0$. Then $M$ admits a Riemannian metric 
of  scalar curvature zero iff $M$ is minimal
and has Kodaira dimension $0$. In particular, 
the Yamabe invariant $Y(M)$ is unachieved 
whenever $\kod (M) =1$. 
\end{prop}

\begin{proof} Any minimal K\"ahler-type
complex surface with 
Kodaira
dimension $0$ has $c_1=0$ in 
real cohomology, and so carries a 
Ricci-flat K\"ahler metric by 
Yau's solution \cite{yau} of the Calabi
conjecture. Such a metric has
$s=0$, and so realizes $Y(M)=0$.

 Conversely \cite{bes}, if a manifold
with $Y(M)=0$ admits a scalar-flat metric,
the metric in question must be Ricci-flat,
since otherwise it could be deformed into
a metric of positive scalar curvature  \cite{bes}. 
In dimension 4, this implies that 
$$2\chi + 3\tau =\frac{1}{4\pi^2}\int_M\left[2|W_+|^2+\frac{s^2}{24}
-\frac{|\stackrel{\circ}{r}|^2}{2}\right]d\mu\geq 0.$$ 
Our complex surface must therefore satisfy $c_1^2\geq 0$. 
If $c_1$ is not a torsion class, the Seiberg-Witten
invariant of $(M,[J])$ is therefore both  metric-independent
and non-zero, even if $b^+(M)=1$.
But the Seiberg-Witten 
estimate \cite{leb2} 
$$\int_Ms^2d\mu \geq 32\pi^2 (c_1^+)^2$$
then forces $c_1$ to be anti-self-dual with 
respect to any scalar-flat metric; thus
$c_1^2\leq 0$, with equality iff $c_1=0$
in real cohomology. Hence $c_1$ is a torsion class,
and our surface must be a minimal surface of
Kodaira dimension $0$. 
\end{proof}

While this makes a satisfactory distinction between
{\em minimal} surfaces of Kodaira dimensions $0$ and $1$,
it is still unsatisfactory for non-minimal surfaces.
However, it is not hard to see that there are
sequences of unit-volume constant-scalar-curvature 
metrics on any blow-up of, say,  a K3 surface 
which bubble off to any given K\"ahler-Einstein 
metric on the minimal model. The metrics in this
sequence have uniformly bounded diameter. 
It seems plausible to conjecture that this does
not happen in the case of  Kodaira dimension $1$:

\begin{conj}
Let $M$ be the underlying 4-manifold of a 
complex algebraic surface of Kodaira dimension 
$\geq 0$. Suppose that there is a sequence
of unit-volume constant-scalar-curvature
metrics on $M$ with uniformly bounded diameter and
$s\nearrow 0$. Then $\kod (M) =0$. 
\end{conj}

Of course, 
beyond simply knowing the sign of $Y(M)$, one would  
like to ask to know its actual value. When $b_1(M)$ is even and 
$\kod (M) \geq 0$, this  question is completely answered by 
Theorems \ref{gtype} and \ref{y0}. For 
algebraic surfaces with  $\kod (M) =-\infty$,
however, the only case in which $Y(M)$ is actually 
known  is for $M=\cp_2$, where $Y(M)$ is actually realized
by the Fubini-Study metric \cite{leb4}. It is, in particular, 
unclear whether blowing up leaves $Y(M)$ unchanged when $\kod (M) 
=-\infty$, although this is what happens when the 
Kodaira dimension is larger.

There are a number of places in this 
paper where we have assumed that the surface in question
has $b_1$ even in order to assure that there is
a non-trivial Seiberg-Witten invariant. However, 
this is superfluous in Kodaira dimension $0$, since 
Kodaira surfaces are symplectic and so have non-trivial
invariants by the work of Taubes \cite{taubes}. 
For Kodaira dimension $1$ the situation is less clear,
but Biquard \cite{biq} has also shown that some non-K\"ahler surfaces
in this class have non-trivial Seiberg-Witten invariants.
Thus it would be surprising if the following 
were not true:

\begin{conj}
Let $M$ be the underlying 4-manifold of a compact 
complex manifold of Kodaira dimension $0$ or $1$.
Then, even if $b_1(M)$ is odd, one still has $Y(M)=0$. 
\end{conj}

Even allowing for the case of $b_1$ odd, the following
is an immediate consequence of the results proven here:

\begin{thm}
Let $(M,J)$ be any compact
complex surface. If $Y(M)  < 0$, then 
$M$ is either  of general type or   of type $VII$. 
\end{thm}

For all known surfaces of type $VII$, 
however, one can check by hand that 
$Y \geq 0$. It therefore seems reasonable to make the following 

\begin{conj}
Let $(M,J)$ be any compact
complex surface. If $Y(M)  < 0$, then 
$M$ is    of general type. 
\end{conj}

Finally, we have accidentally computed the infima of some
other Riemannian functionals by showing that minimal elliptic
surfaces collapse with bounded Ricci curvature.

\begin{prop}
Let $M$ be the underlying 4-manifold of any 
complex elliptic surface. 
Then 
$$\inf_{[g]} \int_M |W_+|^2d\mu =0.$$
Here $W_+$ denotes the self-dual Weyl
tensor, and the infimum is over the set of all 
conformal classes of Riemannian metrics on $M$. 
As a consequence  
$$\inf_{[g]} \int_M |W_-|^2d\mu = -12\pi^2\tau (M).$$
\end{prop}

Indeed, if $M$ is minimal, this follows from
Theorem \ref{rcol} and the Gauss-Bonnet formula 
for 
$(2\chi + 3\tau )(M) = 0$. Since the Burns metric
is anti-self-dual, it is easy to extend this 
to the non-minimal case; cf. \cite{taubes2}. 

By Taubes existence theorem \cite{taubes2} for anti-self-dual
metrics, any complex surface, after being blown up
sufficiently many times, admits metrics with
$W_+=0$. The last result raises, once again, the 
fascinating question of determining 
precisely how many times a given complex surface
must be blown up before this happens. 

\vfill
\noindent
{\bf Acknowledgement.}
The author  warmly
thanks   Michael Anderson, Peter Kronheimer,
and Heberto del Rio
for their helpful comments and suggestions.

\end{document}